\RequirePackage{fix-cm}
\documentclass[smallcondensed]{svjour3}     
\smartqed  
\usepackage{amsfonts}
\usepackage{amssymb}
\usepackage{graphicx}
\usepackage{multirow}
\usepackage[cmex10]{amsmath}
\usepackage{mathtools}

\usepackage{float}
\usepackage{subfig}
\usepackage[utf8]{inputenc}
\usepackage{makeidx}
\usepackage[english]{babel}
\usepackage{verbatim}
\usepackage[autostyle]{csquotes}
\usepackage{authblk}

\usepackage[numbers]{natbib}
\usepackage[bottom]{footmisc}
\usepackage{enumitem}

\graphicspath{{./figures/}}
\begin{document}
	
	\title{Convex Formulation for Planted Quasi-Clique Recovery
	}
	
	
	\author{Sakirudeen A. Abdulsalaam         \and
		M. Montaz Ali 
	}
	
	
	\institute{Sakirudeen A. Abdulsalaam \at
		School of Computer Science and Applied Mathematics, University of the Witwatersrand, Johannesburg, South Africa \\
		\email{sa.abdulsalaam@gmail.com}           
		\and
		M. Montaz Ali \at
		School of Computer Science and Applied Mathematics, University of the Witwatersrand, Johannesburg, South Africa\\
		\email{montaz.ali@wits.ac.za}
	}
	
	\date{Received: date / Accepted: date}

	\maketitle
	
	\begin{abstract}
	In this paper, we consider the planted quasi-clique or $\gamma$-clique problem. This problem is an extension of the well known planted clique problem which is NP-hard. 
	The maximum quasi-clique problem is applicable in community detection, information retrieval and biology. 
	We propose a convex formulation using nuclear norm minimization for planted quasi-clique recovery. We carry out numerical experiments using our convex formulation and the existing mixed integer programming formulations. Results show that the convex formulation performs better than the mixed integer formulations when $\gamma$ is greater than a particular threshold.
	
	\end{abstract}
	
\section{Introduction}
\label{intro}
A clique is the densest subgraph of any undirected graph, $G = (V,E)$. A subgraph, $G[V']$ induced by $V' \subseteq V$, forms a clique if every pair of nodes are adjacent \cite{lucemethod}. The problem of finding the largest clique in a graph is known as the Maximum Clique Problem (MCP) \cite{bomzemaximum,pardalos,pardalosmaximum}.  
The size of the largest clique in $G$ is known as the clique number and it is denoted by $\omega(G)$. Although this problem is NP-hard \cite{garey}, it has been well studied due to its wide applications.  Verily, cliques possess the ideal properties for cohesiveness \cite{pattillocliquerelaxation}. 
However, the requirement that every pair of nodes are adjacent is too confining for some applications \cite{pattillomax}. This motivates the emergence of different clique relaxations. Some of the clique relaxation models emanating from soical network analysis include the $k$-clique, $k$-club, and $k$-plex, see for example, \cite{balasnovel, balascliquerelaxations}. 
A density based relaxation known as quasi-clique or $\gamma$-clique was introduced in \cite{abellomaximum}. Although $\gamma$-clique is the most recent of clique relaxations, it is one of the most popular due to its suitability for a range of applications \cite{veremyevexact}. The quasi-clique model is applicable in community detection \cite{hajekcomputational,pattillocliquerelaxation},
data clustering and data mining \cite{balasnovel,vermanetwork},
information retrieval \cite{terveenconstructing}, protein-protein network \cite{junkeranalysis}, and criminal network analysis \cite{balasthesis,kianianalysis}.

A subgraph, $G[V']$, is a $\gamma$-clique if $|V' \times V' \cap E|/ \binom{|V'|}{2} \geq \gamma$, where $\gamma \in (0, 1]$. 
The problem of finding a $\gamma$-clique with maximum cardinality in $G$ is known as the maximum quasi-clique problem (MQCP). Obviously, the case $\gamma = 1$ is equivalent to the maximum clique problem. This problem has been shown to be NP-hard \cite{pattillomax}. The size of the maximum quasi-clique in $G$ is known as $\gamma$-clique number and we denote it by $\omega_ \gamma(G)$. A number of existing works on finding the maximum $\gamma$-clique focused on developing heuristic methods for finding large quasi-cliques for different instances \cite{abellomassive,bhattacharyyamining,liu2008effective,peimining}. Other various heuristic and enumerative algorithms have recently been developed (see e.g, \cite{marinellilp,miaoellipsoidal,pastukhovmaximum,ribeiroexact,zhouopposition})

The first mathematical model for maximum quasi-clique recovery is the mixed integer programming (MIP) of \cite{pattillomax}. This model has been reformulated in \cite{veremyevexact} to handle larger problems. These are the only known existing maximum quasi-clique recovery models. In this paper, we focus on a special case of this problem, namely, the planted quasi-clique. We propose a novel convex formulation for the planted quasi-clique recovery. We adopt techniques from the matrix decomposition to split the adjacency matrix of the given graph into its low rank and sparse component. We were inspired by the work in \cite{amesnuclear} where a matrix completion strategy was used for planted maximum clique recovery.  Planted clique is a well known problem that has been studied in \cite{alon,amesnuclear,friezenew,kucera}. To the best of our knowledge, this paper presents the first attempt to solve the planted maximum $\gamma$-clique problem. Our numerical experiments show that this approach is more robust than the nuclear norm formulation of \cite{amesnuclear} and  more effective for planted quasi-clique recovery than the mixed integer programming models of \cite{pattillomax,veremyevexact}. 

The rest of this paper is organised as follows. We present the nuclear norm formulation for the planted clique and planted quasi-clique problem in Section \ref{sec:plantedqq}. We briefly present the mixed integer programming formulations for maximum quasi-clique problem in Section \ref{intro2qq}. 
The report of our numerical experiments is presented in Section \ref{sec:num_exp} while the concluding remarks are made in Section \ref{sec:conclusion}.

\section{The planted quasi-clique model}\label{sec:plantedqq}
An instance of the MCP is the planted (hidden) clique problem. 
The problem can be formulated in two different ways namely: the randomized case and the adversarial case. For the randomized case, $n_c$ vertices are chosen at random from $n$ vertices ($n > n_c$) and a clique of size $n_c$ is constructed. The remaining pairs of nodes are then connected depending on a given probability. For the adversarial case, on the contrary, instead of joining the diversionary edges in a probabilistic manner, an adversary is allowed to join the edges. A restriction is placed on the maximum number of edges he can insert so that a clique bigger than the planted clique is not formed. For our planted quasi-clique problem, we have only considered the randomized case. We have formulated the problem using two probablities, namely: $p, p \geq \gamma$, is the probability of an edge existing between two nodes belonging to the planted quasi-clique while $\rho$ is the probability of an edge between the nodes not belonging to the planted quasi-clique. So, summarily, we generate a graph of size $n$ and select $n_c$ nodes randomly ($n_c < n$) and connect them with probability $p$. The remaining $n - n_c$ nodes form the diversionary nodes are they are connected with probability $\rho < p$. The smaller the value of $p$ the more difficulty it is to recover the planted quasi-clique. Conversely, as $\rho$ grows bigger, tending towards $0.5$, the harder it is to recover the planted quasi-clique.

The planted clique problem has previously been studied in \cite{alon,feigefinding,friezenew,kucera}. The following nuclear norm minimization formulation has been recently proposed for solving the planted clique problem \cite{amesnuclear}:
\begin{subequations}\label{mcpnn}
	\begin{align}\label{mcpnn:a}
	&\min ||X||_*, \\
	\label{mcpnn:b}\text{subject to } & \sum_{i \in V} \sum_{j \in V} X_{ij} \geq n_c^2, \\
	\label{mcpnn:c} & X_{ij} = 0 \text{ } \forall \text{ } (i, j) \notin E \text{ and } i \neq j,\\
	\label{mcpnn:d} & X = X^T, \\
	\label{mcpnn:e} &X_{ij} \in [0, 1],	
	\end{align}
\end{subequations}
where $X \in \mathbf{R}^{n \times n}$ and $n_c$ is the size of the planted clique. The nuclear norm is defined as $||X||_* := \sigma_1(X) + \sigma_2(X) + \ldots + \sigma_r(X)$, where $\sigma_i(X), i \in \{1, \ldots, r\}$, are the singlular values and $r$ is the rank of the matrix. We denote model \eqref{mcpnn} as NNM(1) (nuclear norm based model 1).

In our case, we adopt the technique from matrix decomposition \cite{candesrobust} to recover the planted quasi-clique in a graph. The planted quasi-clique problem is a more difficult problem than the planted clique problem, as the latter is a special case of the former. 
The matrix decomposition problem is described as follows. A matrix, $M$, is formed by adding a low rank matrix, $L$, to a sparse matrix, $S$. The objective is to devise a mean to separate $M$ into its low rank and sparse component. Mathematically, we want to solve:
\begin{subequations}\label{mdc}
	\begin{align}
	&\min  rank(L) + ||S||_0,\\
	&\text{subject to } L + S = M,
	\end{align}
\end{subequations}
where $||S||_0 = card(S)$ is the number of non-zero entries of $S$. Both the rank function and $l_0$ minimization are non-convex. However, nuclear norm minimization gives a good approximation of the rank minimization problem \cite{rechtguaranteed}. Furthermore, the matrix $l_1$ norm, defined as $||X||_1 = \sum_{i = 1}^{n_1} \sum_{j = 1}^{n_2}|X_{ij}|$ for $X \in \mathbf{R}^{n_1 \times n_2}$, is a good replacement for the cardinality minimization problem \cite{fornasiernumerical}. Hence, problem \eqref{mdc} can be written as
\begin{subequations}\label{mdc2}
	\begin{align}
	&\min  ||L||_* + ||S||_1,\\
	&\text{subject to } L + S = M.
	\end{align}
\end{subequations}
This problem has applications in facial recognition and image segmentation \cite{candesrobust}. The problem has been studied in \cite{candesrobust, chandrasekaranrank,chenlow}, with application to facial recognition in \cite{candesrobust}. We apply this technique to planted quasi-clique recovery. Our proposed formulation is the following:
\begin{subequations}\label{qcp0}
	\begin{align}\label{qcp0:a}
	&\min ||Q||_* + \lambda ||D||_1 \\
	\label{qcp0:b}\text{subject to } & \sum_i \sum_j Q_{ij} \geq \gamma \eta^2   \\
	\label{qcp0:d}& Q + D = A \\ 
	\label{qcp0:f}& Q_{ij}, D_{ij} \in [0, 1], \quad \eta \in \mathbb{N}, 	
	\end{align}
\end{subequations}
where $Q, D \in \mathbf{R}^{n \times n}$ are matrix variables corresponding to the quasi-clique and the diversionary edges, $A$ is the adjacency matrix of the input graph while the parameter $\gamma \in (0, 1]$ is the desired edge density of the quasi-clique to be recovered. The constraint \eqref{qcp0:b} enusures that the solution satisfies the edge density requirement, while \eqref{qcp0:d} makes sure that the decomposition agrees with the input matrix. $\eta$ is a positive integer value variable that determines the size of the recovered quasi-clique. 

Since we are only interested in $Q$, we can eliminate constraint \eqref{qcp0:d} and write $D = A - Q$. Therefore, \eqref{qcp0} can be reformulated as

\begin{subequations}\label{qcp}
	\begin{align}\label{qcp:a}
	&\min ||Q||_* + \lambda ||A - Q||_1 \\
	\label{qcp:b}\text{subject to } & \sum_i \sum_j Q_{ij} \geq \gamma \eta^2   \\
	\label{qcp:e}& Q_{ij} \in [0, 1], \quad \eta \in \mathbb{N}. 	
	\end{align}
\end{subequations}
We denote model \eqref{qcp} as NNM(5). Following the approach in Appendix A of \cite{chandrasekaranrank}, the semidefinite (SDP) formulation for \eqref{qcp} is the following:
\begin{equation}\label{qcp_sdp}
\begin{aligned}
\text{minimize } &\frac{1}{2}(trace(Z_1) + trace(Z_2)) + \lambda 	\textbf{1}_n^T W \textbf{1}_n,\\
\text{subject to }&  \begin{bmatrix}
Z_1 & Q \\
Q^T &  Z_2
\end{bmatrix} \succeq 0,\\
& -W_{ij} \leq A_{ij} - Q_{ij} \leq W_{ij}, \quad \forall \text{ } (ij),\\
& \sum_i \sum_j Q_{ij} \geq \gamma \eta^2,   \\
& Q_{ij} \in [0, 1], \quad \eta \in \mathbb{N},
\end{aligned}
\end{equation}
where ${\bf 1}_n \in \mathbf{R}^n$ is an $n$-dimensional vector of all entries equal to one and $Z_1, Z_2, W \in \mathbf{R}^{n \times n}$.

Problems \eqref{qcp0} and \eqref{qcp} are convex optimization problems that can be solved using one of the available convex optimization solvers. 
\section{Illustrative Example}
\begin{figure}[H]
\begin{center}
	\includegraphics[width = 0.65 \textwidth]{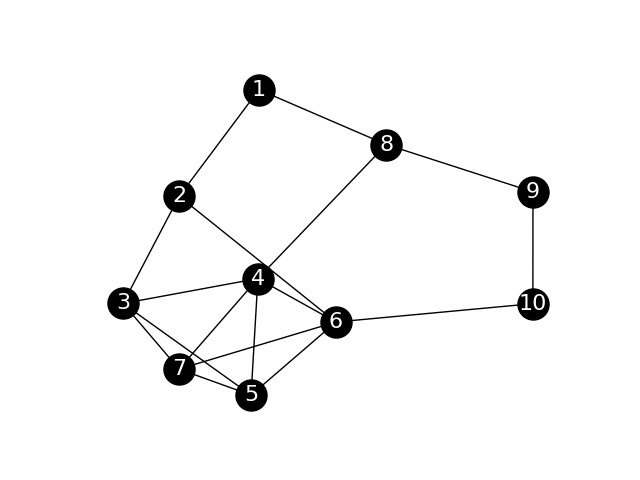}
	\caption{A graph $G$, with the planted quasi-clique using nodes $\{3,4,5,6,7\}$ \label{example1}}
\end{center}
\end{figure}

Suppose the input graph, $G$, containing the planted quasi clique is the graph presented in Figure \ref{example1} and that we want to recover the planted $0.9$-clique from it. $A^*$ is the adjacency matrix of $G$. We add a loop to every node of $G$ to obtain $A$ as the adjacency matrix. This is necessary for the algorithm to be able to recover a low rank submatrix. 
\[
A^* = 
\begin{pmatrix}	
0& 1& 0& 0& 0& 0& 0& 1& 0& 0\\
1& 0& 1& 0& 0& 1& 0& 0& 0& 0\\
0& 1& 0& 1& 1& 0& 1& 0& 0& 0\\
0& 0& 1& 0& 1& 1& 1& 1& 0& 0\\
0& 0& 1& 1& 0& 1& 1& 0& 0& 0\\
0& 1& 0& 1& 1& 0& 1& 0& 0& 1\\
0& 0& 1& 1& 1& 1& 0& 0& 0& 0\\
1& 0& 0& 1& 0& 0& 0& 0& 1& 0\\
0& 0& 0& 0& 0& 0& 0& 1& 0& 1\\
0& 0& 0& 0& 0& 1& 0& 0& 1& 0\\
\end{pmatrix}, \quad A = 
\begin{pmatrix}	
1& 1& 0& 0& 0& 0& 0& 1& 0& 0\\
1& 1& 1& 0& 0& 1& 0& 0& 0& 0\\
0& 1& 1& 1& 1& 0& 1& 0& 0& 0\\
0& 0& 1& 1& 1& 1& 1& 1& 0& 0\\
0& 0& 1& 1& 1& 1& 1& 0& 0& 0\\
0& 1& 0& 1& 1& 1& 1& 0& 0& 1\\
0& 0& 1& 1& 1& 1& 1& 0& 0& 0\\
1& 0& 0& 1& 0& 0& 0& 1& 1& 0\\
0& 0& 0& 0& 0& 0& 0& 1& 1& 1\\
0& 0& 0& 0& 0& 1& 0& 0& 1& 1\\
\end{pmatrix}.
\]
The matrix decomposition algorithm will perform two tasks, namely; completion of the low-rank matrix and separation of the low-rank matrix from the sparse matrix. However, since we are only interested in the low-rank submatrix, we have reformulated the model to suite this purpose. The reformulation has improved the performance of the algorithm in terms of speed. Therefore, for this particular example, we recover $Q$ as the largest rank-one matrix and $\eta = 5$ in this case. This corresponds to the adjacency matrix of the recovered maximum clique. The adjacency matrix of the planted maximum quasi-clique, $Q^*$, can then finally be obtained by setting $Q_{ij} = 0$ if $Q_{ij} \neq A^*_{ij}$.

\[
Q = 
\begin{pmatrix}	
0& 0& 0& 0& 0& 0& 0& 0& 0& 0\\
0& 0& 0& 0& 0& 0& 0& 0& 0& 0\\
0& 0& 1& 1& 1& 1& 1& 0& 0& 0\\
0& 0& 1& 1& 1& 1& 1& 0& 0& 0\\
0& 0& 1& 1& 1& 1& 1& 0& 0& 0\\
0& 0& 1& 1& 1& 1& 1& 0& 0& 0\\
0& 0& 1& 1& 1& 1& 1& 0& 0& 0\\
0& 0& 0& 0& 0& 0& 0& 0& 0& 0\\
0& 0& 0& 0& 0& 0& 0& 0& 0& 0\\
0& 0& 0& 0& 0& 0& 0& 0& 0& 0\\
\end{pmatrix}, \quad Q^* = 
\begin{pmatrix}	
0& 0& 0& 0& 0& 0& 0& 0& 0& 0\\
0& 0& 0& 0& 0& 0& 0& 0& 0& 0\\
0& 0& 0& 1& 1& 0& 1& 0& 0& 0\\
0& 0& 1& 0& 1& 1& 1& 0& 0& 0\\
0& 0& 1& 1& 0& 1& 1& 0& 0& 0\\
0& 0& 0& 1& 1& 0& 1& 0& 0& 0\\
0& 0& 1& 1& 1& 1& 0& 0& 0& 0\\
0& 0& 0& 0& 0& 0& 0& 0& 0& 0\\
0& 0& 0& 0& 0& 0& 0& 0& 0& 0\\
0& 0& 0& 0& 0& 0& 0& 0& 0& 0\\
\end{pmatrix},
\]
\section{Existing Formulations for Maximum Quasi-Clique Problem}\label{intro2qq}
As stated earlier, majority of the existing works on $\gamma$-clique focused on developing heuristics for detection of large quasi-clique. 
The first deterministic solution for $\gamma$-clique recovery problem is based on the linear mixed integer programming model suggested 
in \cite{pattillomax}, where an upper bound was derived. 
For $i \in V$, define $x_i \in \{0, 1\}$ such that $x_i = 1$ if and only if $i \in V'$ and $0$ otherwise, where $V'$ is the vertex set of the maximum quasi-clique. The following linearized MIP formulation was proposed:
\begin{subequations}\label{qcmip2}
	\begin{align}\label{qcmip2:a}
	\omega_\gamma &= 	\max \sum_{i \in V} x_i, \\
	\text{subject to: } & \nonumber \\
	\label{qcmip2:b} \sum_{i \in V} h_ i & \geq  0\\
	\label{qcmip2:c} h_i & \leq  \nu x_i, \quad h_i \geq - \nu x_i \quad \forall \text{ } i \in V,\\
	\label{qcmip2:d}h_i & \geq \gamma x_i + \sum _{j \in V} (A_{ij} - \gamma) x_j  - \nu (1 - x_i) \quad \forall \text{ } i \in V,\\
	\label{qcmip2:e} h_i & \leq \gamma x_i + \sum _{j \in V} (A_{ij} + \gamma) x_j  - \nu (1 - x_i) \quad \forall \text{ } i \in V,\\
	\label{qcmip2:f} x_i &\in \{0, 1\} \quad \forall \text{ } i \in V,
	\end{align}
\end{subequations}
where, $A_{ij}$ are the entries of the adjacency matrix of the graph, $\nu$ is a constant that is large enough and $h_i$ is defined as 
\begin{equation}
h_i = x_i(\gamma x_i + \sum_{(i, j) \in E} (A_{ij} - \gamma)x_j ).
\end{equation} 
\eqref{qcmip2} can only handle problems with small graph size. Because of this drawback, \citeauthor{veremyevexact} \cite{veremyevexact} reformulate this model by defining $z_{ij}$ as a binary variable, such that $z_{ij} = 1$ if and only if $(i,j) \in E \cap (V' \times V')$. In addition, a binary variable $s_t, t = 1, \ldots, |V|$, which determines the size of the quasi-clique is defined. This implies that $s_t = 1$ if and only if $|V'| = t$. With these additional variables and notations, the improved MIP model presented in \cite{veremyevexact} is:
\begin{subequations}\label{qcmip3}
	\begin{align}\label{qcmip3:a}
	&\max \sum_{i \in V} x_i, \\
	\label{qcmip3:b}\text{subject to } & \sum_{(i,j) \in E} z_{ij} \geq \gamma \sum_{t = \omega^l}^{\omega^u} \frac{t(t - 1)}{2}s_t, \\
	\label{qcmip3:c}&z_{ij} \leq x_i, \quad z_{ij} \leq x_j, \text{ } \forall \text{ } (i,j) \in E,\\
	\label{qcmip3:d}& \sum_{i \in V} x_i = \sum_{t = \omega^l}^{\omega^u} t s_t,  \quad \sum_{t = \omega^l}^{\omega^u} s_t = 1, \\ 
	\label{qcmip3:e}& x_i \in \{0, 1\}, z_{ij} \geq 0, \text{ } \forall \text{ } i, j \in V, i < j,\\
	\label{qcmip3:f}& s_t \geq 0, \text{ } \forall \text{ } t \in \{\omega^l, \ldots, \omega^u \},	
	\end{align}
\end{subequations}
where $\omega ^ l$ and $\omega ^ u$ are the upper and lower bound on the size of quasi-clique that could be found in the input graph. These can be set to $0$ and $|V|$ respectively if no estimates are available. The constraint \eqref{qcmip3:b} is the edge density requirement while \eqref{qcmip3:c} ensures that $z_{ij} = 1$ if and only if $i$ and $j$ belong to the quasi-clique. Observe that the left hand side of \eqref{qcmip3:b} can be written as
\begin{equation}\label{6brewrite}
\sum_{(i,j) \in E} z_{ij} = 1/2 \sum_{i \in V} \sum_{j:(i,j) \in E} x_i x_j =  1/2 \sum_{i \in V} \left(x_i \sum_{j:(i,j) \in E} x_j \right).
\end{equation}
Setting $w_i$ to the quantity in the bracket in equation \eqref{6brewrite} above, \eqref{qcmip3} can be reformulated as \cite{veremyevexact}:

\begin{subequations}\label{qcmip4}
	\begin{align}\label{qcmip4:a}
	&\max \sum_{i \in V} x_i \\
	\label{qcmip4:b}\text{subject to } & \sum_{i \in V} w_i \geq \gamma \sum_{t = \omega^l}^{\omega^u} t(t - 1)s_t, \\
	\label{qcmip4:c}&w_i \leq \psi _i x_i, \quad w_i \leq \sum_{j:(i,j) \in E} x_j, \text{ } \forall \text{ } i \in V,\\
	\label{qcmip4:d}& \sum_{i \in V} x_i = \sum_{t = \omega^l}^{\omega^u} t s_t,  \quad \sum_{t = \omega^l}^{\omega^u} s_t = 1, \\ 
	\label{qcmip4:e}& x_i \in \{0, 1\}, z_{ij} \geq 0, \text{ } \forall \text{ } i, j \in V, i < j, ,\\
	\label{qcmip4:f}& s_t \geq 0, \text{ } \forall \text{ } t \in \{\omega^l, \ldots, \omega^u \},	
	\end{align}
\end{subequations}
where $\psi_i$ is a parameter that is sufficiently large. In particular, $\psi_i = deg_G(i)$, where $deg_G(i)$ is the degree of a node $i$ in a given graph, $G$. 

\section{Numerical Experiments}\label{sec:num_exp}
Recall that the planted quasi-clique problem becomes planted clique problem when $\gamma$ is equal to one. We performed numerical experiments with our nuclear norm minimization (NNM) formulation \eqref{qcp} for planted maximum quasi-clique to compare its performance with the existing nuclear norm minimization formulation \eqref{mcpnn} for planted maximum clique recovery. Further, we compare the efficacy of our formulation with the mixed integer programming models \eqref{qcmip2}, \eqref{qcmip3} and \eqref{qcmip4} for quasi cliques. 

\begin{figure}[H]
		\centering
\includegraphics[width = 0.7\textwidth]{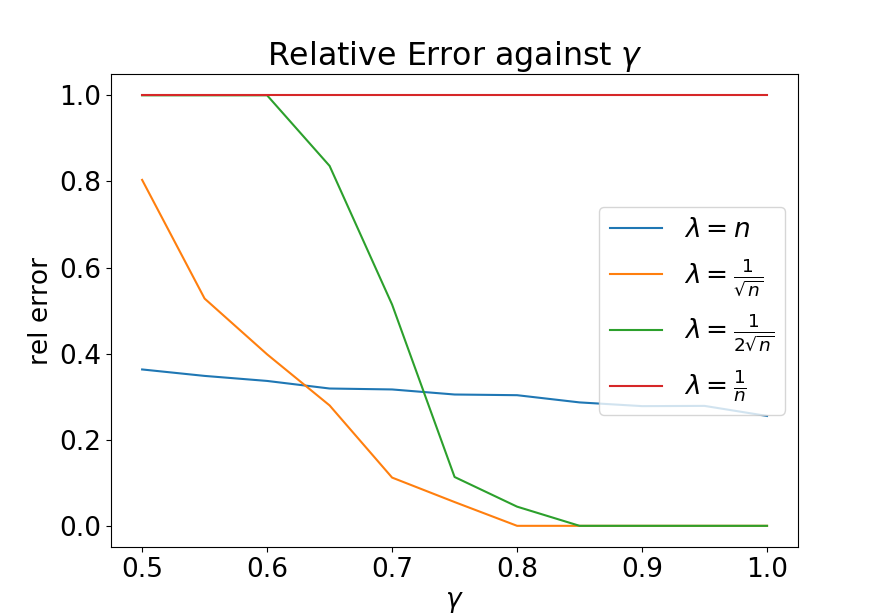}
\caption{Recovery error for different values of $\lambda$.\label{lambdagraph}}
\end{figure}

The experiments have been performed on a HP computer with 16GB Ram and Intel core i7 processor. The machine runs on Debian Linux. The simulations are performed using CVXPY \cite{cvxpy} with NCVX \cite{ncvx}. CVXPY is a python package used to solve convex optimization problems with different solvers, e.g SCS, CVXOPT, and XPRESS. Every instance of the experiment has been carried out ten times and the average result is taken. We have used different values of the regularization parameter, $\lambda$, to ascertain that our choice of $\lambda$ works well for the problem. We have planted a quasi-clique with $\omega_\gamma(G) = 35$ for various values of $\gamma$ in a graph with $50$ nodes. We implement our algorithm for $\lambda = n, \frac{1}{\sqrt{n}}, \frac{1}{2\sqrt{n}}, \frac{1}{n}$. The results is presented in Table \ref{lambdatable} and Figure \ref{lambdagraph}. From Table \ref{lambdatable}, we discover that when $\lambda = n \text{ and } \frac{1}{n}$, the algorithm fails in all instances considered. This finding has been supported by the relative errors in Figure \ref{lambdagraph}, where the relative errors have been calculated using \eqref{error_comp}.

\begin{equation}\label{error_comp}
\text{Relative Error} = \frac{||\text{recovered } \gamma \text{-clique} - \text{planted } \gamma\text{-clique}||_F}{||\text{planted } \gamma\text{-clique}||_F},
\end{equation}
where $F$ denotes the Frobenius norm.
This is contrary to what can be observed when $\lambda = \frac{1}{\sqrt{n}} \text{ and } \frac{1}{2\sqrt{n}}$. For these values, the exact size of planted quasi-clique has been recovered with zero relative error when $\gamma$ approaches $1$ ($\gamma \rightarrow 1$). From this results, we conclude that values of $\frac{1}{\sqrt{n}} \leq \lambda \leq \frac{1}{2\sqrt{n}}$ will work for our model. This is similar to the recommendation in \cite{candesrobust}. This is not surprising since the entries of our matrix are also independent and identically distributed (iid) and hence satisfy the incoherence condition (see \cite{candescomletion_with_noise,candesexactmatrix,candesrobust}). \citeauthor{chandrasekaranrank} \cite{chandrasekaranrank} also contains a heuristic for choosing $\lambda$ and our finding agrees with their result, although our approach is different.  The detailed report of the experiments is as follows.

	\begin{table}[H]
		\centering
	\begin{tabular}{|c|c|c|c|c|c|}\hline
		&Size of the planted quasi-clique& \multicolumn{4}{|c|}{ Size of recovered quasi-clique when}\\\hline
		$\gamma$ &&$\lambda = n$&$\frac{1}{\sqrt{n}}$&$\frac{1}{2\sqrt{n}}$&$\frac{1}{n}$\\\hline
		0.5& 35 &	50&	22.7&	0&	0\\				
		0.55& 35  &  50	&31&	0&	0	\\		
		0.6& 35 	&50	&33.9&	0&	0\\				
		0.65& 35 &	50&	34.5&	12.1&	0	\\			
		0.7& 35 &	50&	34.5&	29.9&	0	\\			
		0.75& 35 &	50&	34.7&	34.6&	0	\\			
		0.8& 35 &	50&	35	&35&	0		\\		
		0.85& 35 &	50&	35&	35&	0\\
		0.9&35&		50	&	35	&	35&	0\\
		0.95&35	&	50	&	35	&	35	&	0\\
		1&	35&	50&		35&		35&		0\\\hline
	\end{tabular}
\caption{Quasi-clique recovery for different values of $\lambda$.\label{lambdatable}}
\end{table}

\subsection{Comparison between NNM\eqref{mcpnn}  and NNM\eqref{qcp}}
\begin{figure}[htbp]
	\centering
	\subfloat[Probability of planted clique recovery with clique size $ = 50$ \label{pr50}]{
		\includegraphics[width = 0.45 \textwidth]{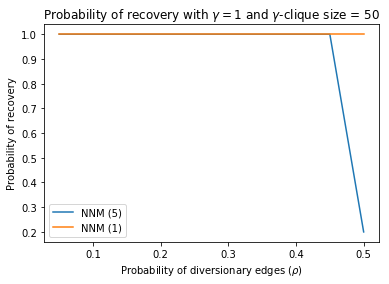}}
		\hfill
	\subfloat[Probability of planted clique recovery with clique size $= 80$ \label{pr80}]{
		\includegraphics[width =0.45 \textwidth]{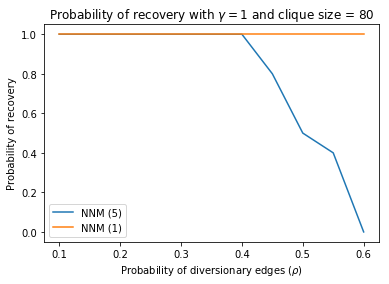}}
	\caption{Planted clique recovery from a graph with 100 nodes with varied probability of adding a diversionary edge ($\rho$).}\label{clicpr}
\end{figure}

\begin{figure}[htbp]
	\centering
	\subfloat[CPU time comparison with clique size = 50 \label{cliccput50}]{
		\includegraphics[width = 0.45 \textwidth]{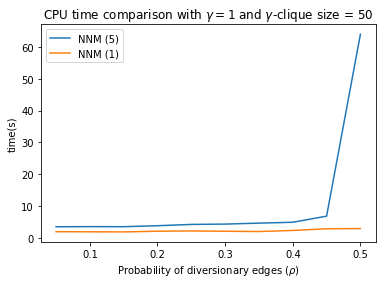}}
	\hfill
	\subfloat[CPU time comparison with clique size = 80 \label{cliccput80}]{
		\includegraphics[width =0.45 \textwidth]{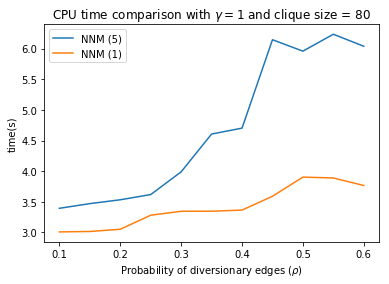}}
		\caption{CPU time comparison for planted clique recovery from a graph with 100 nodes with varied number of diversionary edges.}\label{clictime}
\end{figure}

Our model, NNM\eqref{qcp}, represents the planted maximum clique model when $\gamma = 1$. Hence, we compare its performance with NNM\eqref{mcpnn}.
We have considered a graph with 100 nodes in this experiment. We have planted clique of size $80$ and $50$ and then varied the probability, $\rho$, of an edge existing between the remaining nodes. The results of this experiment is contained in Figure \ref{clicpr}. In both cases considered, the results show that both \eqref{mcpnn} and \eqref{qcp} recover planted clique perfectly when the probability of adding a diversionary edge is below certain threshold (roughly $0.45$). However, \eqref{qcp} fails to perfectly recover the planted clique when this threshold is exceeded while \eqref{mcpnn} still solves the problem perfectly. We  have observed that the presence of constraint \eqref{mcpnn:c} enables finding the largest rank-one submatrix in the input matrix easier in the formulation \eqref{mcpnn}. However, this constraint can not be imposed in the case of \eqref{qcp}, otherwise, solving planted quasi-clique problem with the formulation will be impossible. Figure \ref{clictime} represents the CPU times for NNM\eqref{mcpnn}  and NNM\eqref{qcp} with planted clique of size $\omega_\gamma(G) = 50 \text{ and } 80$ and $\gamma = 1$. It can be observed from the figures that NNM\eqref{mcpnn} is more efficient than NNM\eqref{qcp} in this case.

\begin{figure}[htbp]
	\centering
	\subfloat[Recovery probability of planted quasi-clique with $\gamma = 0.99$ \label{qclicpr80}]{
		\includegraphics[width = 0.45 \textwidth]{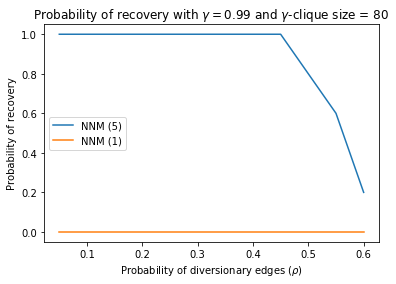}}
	\hfill
	\subfloat[CPU time comparison for planted quasi-clique \label{qcliccput80}]{
		\includegraphics[width =0.45 \textwidth]{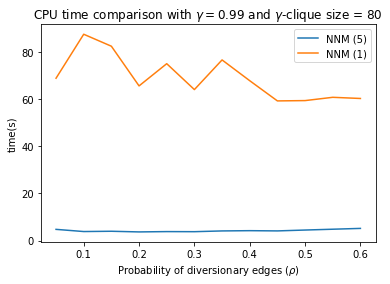}}
	\caption{Comparison of the recovery probability and CPU time of NNM for planted quasi-clique.} \label{qclic80}
\end{figure}

Figure \ref{qclic80} shows the performance of NNM\eqref{mcpnn} compared with NNM\eqref{qcp} in finding quasi-clique (with $\gamma < 1$). We have observed from Figure \ref{qclicpr80} that despite the fact that the planted quasi-clique that we have considered for this case has very few missing edges ($\gamma = 0.99$), NNM\eqref{mcpnn} failed to recover the quasi-clique for every trial. However, NNM\eqref{qcp} produced similar result as the case $\gamma = 1$ (see Figure \ref{pr80} and \ref{qclicpr80}).  In addition, NNM\eqref{qcp} is, by far, more efficient than  NNM\eqref{mcpnn} for the case $\gamma = 0.99$ (see Figure \ref{qcliccput80}).

\subsection{Maximum  Quasi-Clique Recovery}
\subsubsection{The planted case}
Two types of experiment have been performed in this case. In the first case, we have checked whether the recovered quasi-clique satisfies the edge density requirement or not. The second experiment focuses on the size of the recovered quasi-clique, i.e, to examine whether the size of the planted quasi-clique ($n_c$) is the same as the size of the recovered quasi-clique ($\eta$). The detailed report of both experiments is as follows.

The goal of the first experiment is to examine the error in the edge density of the recovered $\gamma$-clique with respect to the edge density of the planted maximum $\gamma$-clique. We have computed the relative error between the edge density of the recovered $\gamma$-clique and the edge density of the planted $\gamma$-clique (i.e, the expected edge density) for various $\gamma$. All the errors computed in this section are relative errors.

\begin{figure}[htbp]
	\centering
	\subfloat[CPU time for planted quasi-clique recovery for a graph with 50 nodes \label{cput50}]{
		\includegraphics[width = 0.45 \textwidth]{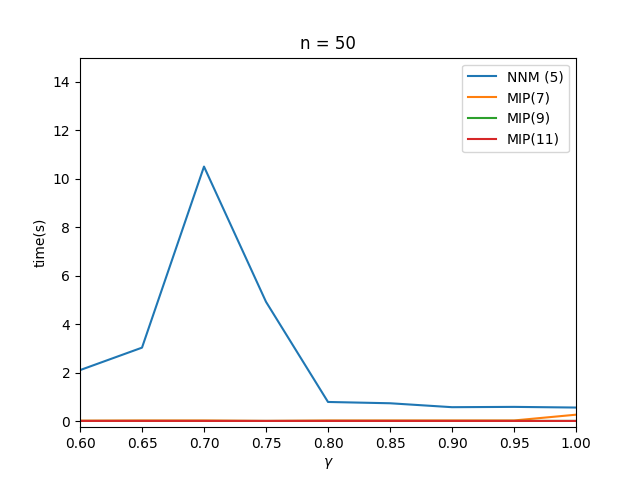}}
	\hfill
	\subfloat[CPU time for planted quasi-clique recovery for a graph with 100 nodes 
\label{cput100}]{
		\includegraphics[width =0.45 \textwidth]{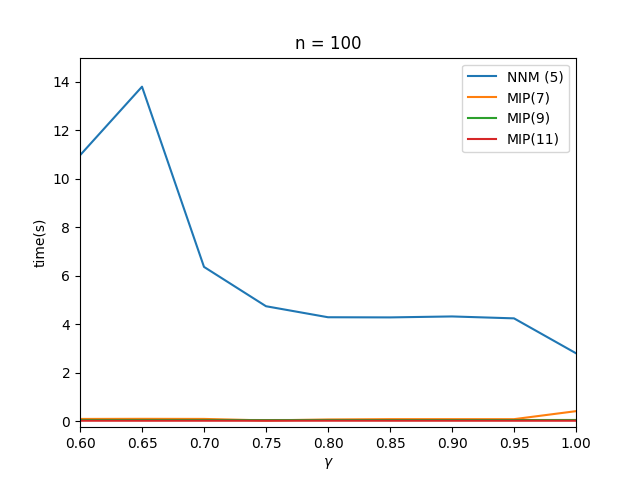}}
		\caption{Comparison of the CPU time for the MIP and NNM methods}
	\label{cput}		
\end{figure}

We have considered again graphs with $50$ and $100$ nodes for this case with planted $\gamma$-cliques of sizes $40$ and $80$, respectively. The planted $\gamma$-clique corresponds to a dense submatrix of the $50 \times 50$ (respectively, $100 \times 100$) input matrix with $40$ (respectively, $80$) nonzero rows/columns. We have varied the edge density of the planted $\gamma$-clique by setting $p = 0.6, 0.65, 0.7, \ldots, 1$. The probability, $p$, determines whether an edge will exist between two nodes in the planted quasi-clique. The smaller the $p$, the fewer the edges and consequently, the more difficult it is to recover what is planted. The setup follows the Stochastic Block Model (SBM) \cite{leereview}. Detail is as follows. For the case $n = 50$, we generate a $50 \times 50$ symmetric matrix, $M$, with zero entries. We choose a $40 \times 40$ submatrix of this matrix and assign $1$ to its indices with probability $0.6$ (suppose $p = 0.6$), using Bernoulli trial. This forms the dense component of the input matrix (the planted $\gamma$-clique). The entries of the remaining $10$ rows and columns are also assigned values $1$ but with a much smaller probability, (say $\rho = 0.2$). This forms the sparse component of the matrix (or the random noise). The goal is to recover the dense submatrix from the input matrix. The results of these experiments are reported in Tables \ref{table50} and \ref{table100}. In both Tables, columns $2 - 4$ contain errors in edge density of the planted quasi-cliques recovered using the MIP models \eqref{qcmip2}, \eqref{qcmip3} and \eqref{qcmip4} while column $5$ contains the errors in edge density of the $\gamma$-clique recovered using our nuclear norm minimization approach, NNM\eqref{qcp}. The relative error here shows the disparity in the densities of what is planted and what is recovered. If these edge densities coincide, i.e, if the edge densities of what is planted and the recovered quasi-clique are equal, the relative error with respect to the Frobenius norm will be zero.  When $n = 50$, MIP\eqref{qcmip2} performed better than the two other MIP models for all values of $\gamma$. However, our model NNM(5) has exhibited the best performance when $\gamma \geq 0.75$. For graphs with $100$ nodes (see Table \ref{table100}), MIP\eqref{qcmip4} performed better than other MIP models except for when $\gamma$ is equal to $0.75, 0.95$ and $1$ where MIP\eqref{qcmip2} has shown better performances. Nevertheless, when $\gamma \geq 0.7$, NNM(5) outperformed all the mixed integer programs. One can also infer from Tables \ref{table50} and \ref{table100} that as the graph size increases, the lower bound on $\gamma$ for perfect recovery by NNM\eqref{qcp} decreases. Figure \ref{cput} shows the CPU time for each of the methods for the experiments reported in Tables \ref{table50} and \ref{table100}. Our off-the-shelf solver, splitting conic solver (SCS) \cite{scs}, is faster than the popular SDP solvers like SeDumi \cite{sturmusing} and SDP3 \cite{tohsdpt3}. However, it is not as efficient as the well-developed FICO XPRESS optimizer used to solve the MIP models. Nonethless, as $\gamma$ increases, there is a drastic drop in the CPU time for NNM\eqref{qcp} in both instances.

\begin{table}[!htbp]
	\begin{center}
		\caption{Errors in the edge density of the planted maximum $\gamma$-clique recovery for a graph with $50$ nodes}\label{table50}
		\begin{tabular}{*{5}{|p{2cm}}|}
			\hline
			\multirow{2}{*}{\bfseries $\gamma$}	&\multicolumn{4}{|c|}{\bf Recovery Error}\\
			\cline{2-5}
			& \bfseries \bf MIP\eqref{qcmip2} &\bf MIP\eqref{qcmip3} &\bf MIP\eqref{qcmip4} &\bf NNM\eqref{qcp}\\ \hline
			0.6	&	0.0922&	0.1279&	0.2093&	0.3085\\	\hline
			0.65	&	0.0688&	0.1607&	0.1897&	0.2053\\	\hline
			0.7	&	0.0809&	0.1646&	0.2086&	0.1170\\	\hline
			0.75	&	0.0809&	0.1935&	0.1558&	0.0230\\	\hline
			0.8	&	0.087&	0.2044&	0.1526&	0\\	\hline
			0.85	&	0.02&	0.2265&	0.1299&	0\\	\hline
			0.9	&	0.0215&	0.2233&	0.1806&	0\\	\hline
			0.95	&	0&	0.2245&	0.2026&	0\\	\hline
			1	&	0&	0.2236&	0.2434&	0\\	\hline
		\end{tabular}
	\end{center}
\end{table}
\begin{table}[!htbp]
	\begin{center}
		\caption{Errors in the edge density of the planted maximum $\gamma$-clique recovery for a graph with $100$ nodes}\label{table100}
		\begin{tabular}{*{5}{|p{2cm}}|}
			\hline
			\multirow{2}{*}{\bfseries $\gamma$}	&\multicolumn{4}{|c|}{\bf Recovery Error}\\
			\cline{2-5}
			&  \bfseries \bf MIP\eqref{qcmip2} &\bf MIP\eqref{qcmip3} &\bf MIP\eqref{qcmip4} &\bf NNM\eqref{qcp}\\ \hline
			0.6	&	0.0916&	0.0736&	0.054&	0.2424\\	\hline
			0.65	&	0.0892&	0.0843&	0.0492&	0.1012\\	\hline
			0.7	&	0.0879&	0.0719&	0.0634&	0.0131\\	\hline
			0.75	&	0.0879&	0.1378&	0.0905&	0\\	\hline
			0.8	&	0.0829&	0.1001&	0.0766&	0\\	\hline
			0.85	&	0.0783&	0.1563&	0.0603&	0\\	\hline
			0.9	&	0.0817&	0.1144&	0.0975&	0\\	\hline
			0.95	&	0.0735&	0.1432&	0.1376&	0\\	\hline
			1	&	0.0694&	0.1581&	0.1717&	0\\	\hline
		\end{tabular}
	\end{center}
\end{table}

The second experiment was to find out if the number of nodes in the planted quasi-cliques, $n_c$, is the same as the number of nodes in the recovered quasi-cliques, $\eta$. For this experiment, we have considered graphs of sizes $n = 50, 100, \ldots, 250$ and $\gamma = 0.6, 0.7, \ldots, 1$. We chose the size of the planted quasi clique, $n_c$, to be  $0.8 \times n$. We have again run the experiment $10$ times for each case and averaged the recovered quasi-clique size. The results obtained are presented in Table \ref{tablequasicliquesizeerror}. The first column under each method contains the average size of recovered quasi-clique using the method while the second contains the relative error for the method. We compute relative error in this case using 
$$\frac{|\text{size of the recovered quasi-clique} - \text{size of the planted quasi-clique}|}{|\text{size of the planted quasi-clique}|}.$$ 
Clearly, if the size of the recovered quasi-clique is equal to the size of the planted quasi-clique, this error will be equal to zero. As shown in the last column of Table \ref{tablequasicliquesizeerror}, the relative errors in the size of quasi-clique recovered via NMM\eqref{qcp} are all zero since $n_c = \eta$ throughout. This shows that the convex formulation always returns correct planted quasi-clique size. MIP\eqref{qcmip2} has the overall worst performance in this experiment. Based on the results in Table \ref{table50} and \ref{tablequasicliquesizeerror}, when $\gamma > 0.75$, $n_c = \eta$ and the error in edge density is equal to zero. 
This implies that our convex formulation perfectly recovers maximum planted quasi-clique when $\gamma > 0.75$ for $n \geq 50$ and $n_c$ large enough. 

\begin{table}[!htbp]
	\begin{center}
		\caption{Errors in the size of planted maximum $\gamma$-clique recovered using different methods for $\gamma$ ranging from $0.6$ to $1$. $n$ is the graph size while $n_c$ is the size of the planted $\gamma$-clique.}\label{tablequasicliquesizeerror}
		\begin{tabular}{*{10}{|c}|}
			\hline
			\multicolumn{2}{|c|}{$\bf \gamma = 0.6$}	&\multicolumn{8}{|c|}{\bf Average Recovered Quasi-clique size/Relative Error}\\
			\hline
		$\bf n$	& $\bf n_c$ &\multicolumn{2}{|c|}{ \bf MIP\eqref{qcmip2}} &\multicolumn{2}{|c|}{\bf MIP\eqref{qcmip3}} & \multicolumn{2}{|c|}{\bf MIP\eqref{qcmip4}} &\multicolumn{2}{|c|}{\bf NNM\eqref{qcp}}\\\hline
			&&&&&&&&$\eta$& \\ \hline
			50	&40		& 41		&0.025		&40.4	&0.01	&40.8	&0.02	&40	&0 \\ \hline
			100	&80		&81.7	&0.021	&80.8	&0.01	&80.8	&0.01	&80		&0 \\ \hline
			150	&120	&122.7	&0.023	&121	&0.008	&120.6	&0.005	&120	&0 \\ \hline
			200	&160	&163.4&	0.021	&161.6	&0.01	&160.9	&0.006	&160	&0 \\ \hline
			250	&200	&204.4	&0.022	&201.9	&0.01	&200.7	&0.003	&200	&0 \\ \hline
			\multicolumn{2}{|c|}{$\bf \gamma = 0.7$}	&\multicolumn{8}{|c|}{}\\	\hline
			50	&40		&40.5	&0.013	&40.1	&0.003	&40.4	&0.01	&40	&0 \\ \hline
			100	&80		&81.3	&0.016	&80.4	&0.005	&80.5	&0.006	&80	&0\\ \hline
			150	&120	&122.4	&0.02	&121	&0.008	&120.7	&0.006	&120	&	0\\ \hline
			200	&160	&163	&0.019	&161	&0.006	&160.6	&0.004	&160	&	0\\ \hline
			250	&200	&204	&0.02	&201.5	&0.008	&200.4	&0.002	&200	&	0\\ \hline
			\multicolumn{2}{|c|}{$\bf \gamma = 0.8$}	&\multicolumn{8}{|c|}{}\\	\hline
			50	&40		&40.2	&0.005	&40		&0		&40.5	&0.013	&40	&0 \\ \hline
			100	&80		&81	&0.025	&80.2	&0.003	&80.2	&0.003	&80		&0\\ \hline
			150	&120	&123	&0.025	&120.5	&0.004	&120.3	&0.002	&120	&0\\ \hline
			200	&160	&163.2	&0.02	&161.2	&0.007	&160.4	&0.003	&160	&0\\ \hline
			250	&200	& 204	& 0.02	&201.3	&0.007	&200.4	&0.002	&200	&0\\ \hline
			\multicolumn{2}{|c|}{$\bf \gamma = 0.9$}	&\multicolumn{8}{|c|}{}\\	\hline
			50	&40		&40.6	&0.015	&40		&0		&40.5	&0.013	&40		&0\\ \hline
			100	&80		&81.1	&0.014	&80.2	&0.003	&80.2	&0.003	&80		&0\\ \hline
			150	&120	&122	&0.017	&120.3	&0.002	&120.1	&0.001	&120	&0\\ \hline
			200	&160	&163	&0.019	&160.9	&0.006	&160.3	&0.002	&160	&0\\ \hline
			250	&200	&203.9	&0.02	&200.9	&0.005	&200.1	&0		&200	&0\\ \hline
			\multicolumn{2}{|c|}{$\bf \gamma = 1$}	&\multicolumn{8}{|c|}{}\\	\hline
			50	&40		&40		&0		&40		&0		&40		&0	&40		&0 \\ \hline
			100	&80		&81		&0.013	&80		&0		&80		&0	&80		&0\\ \hline
			150	&120	&122	&0.017	&120.3	&0.002	&120	&0	&120	&0\\ \hline
			200	&160	&163	&0.019	&160.9	&0.006	&160	&0	&160	&0\\ \hline
			250	&200	&204	&0.02	&201	&0.005	&200	&0	&200	&0\\ \hline
			
		\end{tabular}
	\end{center}
\end{table}
\subsubsection{Recovery from random graphs}
\begin{figure}[H]
	\centering
	\subfloat[CPU time for power law graph with 50 nodes \label{powerlawcput50}]{
		\includegraphics[width = 0.45 \textwidth]{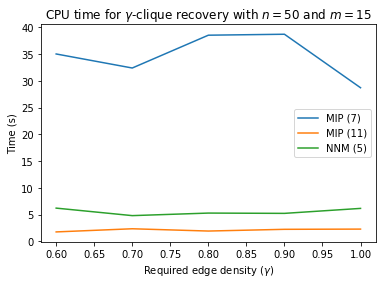}}
		\hfill
	\subfloat[CPU time for power law graph with 100 nodes \label{powerlawcput100}]{
		\includegraphics[width = 0.45 \textwidth]{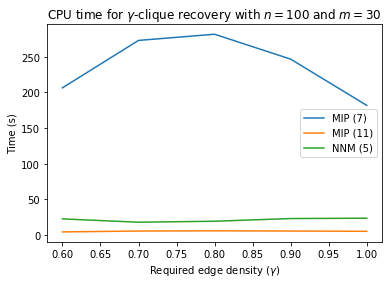}}
	\caption{CPU time comparison for quasi-clique recovery from random graph}\label{powerlawcput}
\end{figure}

Our last experiment focuses on checking the performance of our model in a scenario that mirrors real-life situation. It has been observed that real networks obey some scaling laws rather than being completely random. Hence, the well-known Erdos Renyi random graph, where edges are generated with a constant probability with degree distribution following a Poisson law, may not be suitable. Hence, we have generated our random graph using the preferential attachment model of Barabasi-Albert \cite{barabasiemergence}. The degree distribution of these graphs follow power-law. In this setting, the rate, $\Pi (k)$, with which a node with $k$ edges acquires new edges is a monotonically increasing function of $k$. The time evolution of the degree $k_i$ of node $i$ can be obtained from the first-order ordinary differential equation \cite{jeongmeasuring}:
\begin{equation}
\dfrac{d k _i}{dt} = m \Pi (k_i),
\end{equation} 
where $m$ is a constant; it is the number of edges to attach from a new node to the existing nodes. We have considered graphs with $50$ and $100$ nodes with $m$ set to $15$ and $30$, respectively. The results of this experiment are presented in Table 4. From Table 4(a) and 4(b), it can be observed that MIP\eqref{qcmip2} returns the largest quasi-clique while our NNM\eqref{qcp} returns quasi-cliques with the smallest size. However, for $\gamma \geq 0.8$, our formulation and MIP\eqref{qcmip4} return similar results. Recall, from the first experiment of Section 4.2.1, that the recovery error of our formulation is zero for $\gamma \geq 0.8$. Unfortunately, since the quasi-cliques in this case have not been planted, computing the error in the recovered quasi-clique is not straight-forward. MIP\eqref{qcmip2} has the worst performance in terms of CPU time for this experiment while MIP\eqref{qcmip4} has the best performance of the three formulations compared (see Figure \ref{powerlawcput}). Also, both MIP\eqref{qcmip4} and NNM\eqref{qcp} show no significance difference in CPU time for various value of $\gamma$.

\begin{table}[!htbp]
		\caption{Quasi-cliqe recovery from random graph}
	\subfloat[Quasi-clique recovery from a powerlaw graph with $n = 50$ and $m = 15$.]{
		\centering
		\begin{tabular}{|c|c|c|c|}\hline
			$\gamma$ & MIP\eqref{qcmip2} & MIP\eqref{qcmip4} & NNM\eqref{qcp} \\ \hline
			0.6&	38&	39&	35\\
			0.7&	34&	35&	32\\
			0.8&	32&	30&	29\\
			0.9&	30&	29&	28\\
			1&		29 &	27&	27\\\hline
		\end{tabular}}
		\hfill
		\subfloat[Quasi-clique recovery from a powerlaw graph with $n = 100$ and $m = 30$.]{
		\centering
		\begin{tabular}{|c|c|c|c|}\hline
			$\gamma$ & MIP\eqref{qcmip2} & MIP\eqref{qcmip4} & NNM\eqref{qcp} \\ \hline
			0.6&	76&	78&	68\\
			0.7&	68&	69&	65\\	
			0.8&	64&	63&	64\\	
			0.9	&  61&	58&	59\\	
			1&	58&	53& 54\\\hline
		\end{tabular}}

\end{table}

\section{Conclusion}\label{sec:conclusion}
We have studied the planted quasi-clique problem in this paper. 
We have considered a matrix decomposition type of mathematical formulation for the problem. We have used this formulation to solve the planted maximum quasi-clique problem. We have shown, experimentally, the range of values of the regularization parameter, $\lambda$, that works for the model. We have numerically established the superiority of our formulation over the nuclear norm minimization model in \cite{amesnuclear} by solving a wider range the problem and the three existing mixed integer programming formulations in terms of effectiveness. 
Our future research will be to establish the theoretical guarantee for perfect recovery and providing a bound on $\gamma$ for which recovery is guaranteed. There are some special algorithms developed for nuclear norm minimization and low-rank plus sparse matrix recovery like the iterative singular value thresholding \cite{caisingular}, accelerated proximal gradient \cite{tohaccelerated} and the alternating direction method \cite{yuansparse,ganeshfast}. It will be interesting to implement these algorithms for planted quasi-clique recovery to compare their performances with the SCS used for this work. Lastly, there is no theory to explain why the values of $\lambda$ that work do. It will be interesting have a better understanding of why they do.


	%
	%

	
\bibliographystyle{plainnat}
\bibliography{reference.bib}	

\begin{thebibliography}{50}
\providecommand{\natexlab}[1]{#1}
\providecommand{\url}[1]{\texttt{#1}}
\expandafter\ifx\csname urlstyle\endcsname\relax
  \providecommand{\doi}[1]{doi: #1}\else
  \providecommand{\doi}{doi: \begingroup \urlstyle{rm}\Url}\fi

\bibitem[Abello et~al.(1999)Abello, Pardalos, and Resende]{abellomaximum}
James Abello, PM~Pardalos, and MGC Resende.
\newblock On maximum clique problems in very large graphs.
\newblock \emph{DIMACS series}, 50:\penalty0 119--130, 1999.

\bibitem[Abello et~al.(2002)Abello, Resende, and Sudarsky]{abellomassive}
James Abello, Mauricio~GC Resende, and Sandra Sudarsky.
\newblock Massive quasi-clique detection.
\newblock In \emph{Latin American Symposium on Theoretical Informatics}, pages
  598--612. Springer, 2002.

\bibitem[Alon et~al.(1998)Alon, Krivelevich, and Sudakov]{alon}
Noga Alon, Michael Krivelevich, and Benny Sudakov.
\newblock Finding a large hidden clique in a random graph.
\newblock \emph{Random Structures and Algorithms}, 13\penalty0 (3-4):\penalty0
  457--466, 1998.

\bibitem[Ames and Vavasis(2011)]{amesnuclear}
Brendan~PW Ames and Stephen~A Vavasis.
\newblock Nuclear norm minimization for the planted clique and biclique
  problems.
\newblock \emph{Mathematical programming}, 129\penalty0 (1):\penalty0 69--89,
  2011.

\bibitem[Balabhaskar(2007)]{balasthesis}
Balasundaram Balabhaskar.
\newblock \emph{Graph theoretic generalization of clique: Optimization and
  extensions}.
\newblock PhD thesis, PhD thesis, Texas A \& M University, 2007.

\bibitem[Balasundaram et~al.(2005)Balasundaram, Butenko, and
  Trukhanov]{balasnovel}
Balabhaskar Balasundaram, Sergiy Butenko, and Svyatoslav Trukhanov.
\newblock Novel approaches for analyzing biological networks.
\newblock \emph{Journal of Combinatorial Optimization}, 10\penalty0
  (1):\penalty0 23--39, 2005.

\bibitem[Balasundaram et~al.(2011)Balasundaram, Butenko, and
  Hicks]{balascliquerelaxations}
Balabhaskar Balasundaram, Sergiy Butenko, and Illya~V Hicks.
\newblock Clique relaxations in social network analysis: The maximum $k$-plex
  problem.
\newblock \emph{Operations Research}, 59\penalty0 (1):\penalty0 133--142, 2011.

\bibitem[Barab{\'a}si and Albert(1999)]{barabasiemergence}
Albert-L{\'a}szl{\'o} Barab{\'a}si and R{\'e}ka Albert.
\newblock Emergence of scaling in random networks.
\newblock \emph{science}, 286\penalty0 (5439):\penalty0 509--512, 1999.

\bibitem[Bhattacharyya and Bandyopadhyay(2009)]{bhattacharyyamining}
Malay Bhattacharyya and Sanghamitra Bandyopadhyay.
\newblock Mining the largest quasi-clique in human protein interactome.
\newblock In \emph{2009 International Conference on Adaptive and Intelligent
  Systems}, pages 194--199. IEEE, 2009.

\bibitem[Bomze et~al.(1999)Bomze, Budinich, Pardalos, and
  Pelillo]{bomzemaximum}
Immanuel~M Bomze, Marco Budinich, Panos~M Pardalos, and Marcello Pelillo.
\newblock The maximum clique problem.
\newblock In \emph{Handbook of combinatorial optimization}, pages 1--74.
  Springer, 1999.

\bibitem[Cai et~al.(2010)Cai, Cand{\`e}s, and Shen]{caisingular}
Jian-Feng Cai, Emmanuel~J Cand{\`e}s, and Zuowei Shen.
\newblock A singular value thresholding algorithm for matrix completion.
\newblock \emph{SIAM Journal on optimization}, 20\penalty0 (4):\penalty0
  1956--1982, 2010.

\bibitem[Candes and Plan(2010)]{candescomletion_with_noise}
Emmanuel~J Candes and Yaniv Plan.
\newblock Matrix completion with noise.
\newblock \emph{Proceedings of the IEEE}, 98\penalty0 (6):\penalty0 925--936,
  2010.

\bibitem[Cand{\'e}s and Recht(2009)]{candesexactmatrix}
Emmanuel~J Cand{\'e}s and Benjamin Recht.
\newblock Exact matrix completion via convex optimization.
\newblock \emph{Foundations of Computational mathematics}, 9\penalty0
  (6):\penalty0 717--772, 2009.

\bibitem[Cand{\`e}s et~al.(2011)Cand{\`e}s, Li, Ma, and Wright]{candesrobust}
Emmanuel~J Cand{\`e}s, Xiaodong Li, Yi~Ma, and John Wright.
\newblock Robust principal component analysis?
\newblock \emph{Journal of the ACM (JACM)}, 58\penalty0 (3):\penalty0 11:1 --
  11:37, 2011.

\bibitem[Chandrasekaran et~al.(2011)Chandrasekaran, Sanghavi, Parrilo, and
  Willsky]{chandrasekaranrank}
Venkat Chandrasekaran, Sujay Sanghavi, Pablo~A Parrilo, and Alan~S Willsky.
\newblock Rank-sparsity incoherence for matrix decomposition.
\newblock \emph{SIAM Journal on Optimization}, 21\penalty0 (2):\penalty0
  572--596, 2011.

\bibitem[Chen et~al.(2013)Chen, Jalali, Sanghavi, and Caramanis]{chenlow}
Yudong Chen, Ali Jalali, Sujay Sanghavi, and Constantine Caramanis.
\newblock Low-rank matrix recovery from errors and erasures.
\newblock \emph{IEEE Transactions on Information Theory}, 59\penalty0
  (7):\penalty0 4324--4337, 2013.

\bibitem[Diamond and Boyd(2016)]{cvxpy}
Steven Diamond and Stephen Boyd.
\newblock {CVXPY}: A {P}ython-embedded modeling language for convex
  optimization.
\newblock \emph{Journal of Machine Learning Research}, 17\penalty0
  (83):\penalty0 1--5, 2016.

\bibitem[Diamond et~al.(2018)Diamond, Takapoui, and Boyd]{ncvx}
Steven Diamond, Reza Takapoui, and S~Boyd.
\newblock A general system for heuristic minimization of convex functions over
  non-convex sets.
\newblock \emph{Optimization Methods and Software}, 33\penalty0 (1):\penalty0
  165--193, 2018.

\bibitem[Feige and Krauthgamer(2000)]{feigefinding}
Uriel Feige and Robert Krauthgamer.
\newblock Finding and certifying a large hidden clique in a semirandom graph.
\newblock \emph{Random Structures \& Algorithms}, 16\penalty0 (2):\penalty0
  195--208, 2000.

\bibitem[Fornasier(2010)]{fornasiernumerical}
Massimo Fornasier.
\newblock Numerical methods for sparse recovery.
\newblock \emph{Theoretical foundations and numerical methods for sparse
  recovery}, 14:\penalty0 93--200, 2010.

\bibitem[Frieze and Kannan(2008)]{friezenew}
Alan Frieze and Ravi Kannan.
\newblock {A new approach to the planted clique problem}.
\newblock In \emph{IARCS Annual Conference on Foundations of Software
  Technology and Theoretical Computer Science}, volume~2 of \emph{Leibniz
  International Proceedings in Informatics (LIPIcs)}, pages 187--198, 2008.
\newblock ISBN 978-3-939897-08-8.

\bibitem[Ganesh et~al.(2009)Ganesh, Lin, Wright, Wu, Chen, and Ma]{ganeshfast}
Arvind Ganesh, Zhouchen Lin, John Wright, Leqin Wu, Minming Chen, and Yi~Ma.
\newblock Fast algorithms for recovering a corrupted low-rank matrix.
\newblock In \emph{2009 3rd IEEE International Workshop on Computational
  Advances in Multi-Sensor Adaptive Processing (CAMSAP)}, pages 213--216. IEEE,
  2009.

\bibitem[Garey and Johnson(2002)]{garey}
Michael~R Garey and David~S Johnson.
\newblock \emph{Computers and intractability}, volume~29.
\newblock wh freeman New York, 2002.

\bibitem[Hajek et~al.(2015)Hajek, Wu, and Xu]{hajekcomputational}
Bruce Hajek, Yihong Wu, and Jiaming Xu.
\newblock Computational lower bounds for community detection on random graphs.
\newblock In \emph{Conference on Learning Theory}, pages 899--928, 2015.

\bibitem[Jeong et~al.(2003)Jeong, N{\'{e}}da, and
  Barab{\'{a}}si]{jeongmeasuring}
H~Jeong, Z~N{\'{e}}da, and A.~L Barab{\'{a}}si.
\newblock Measuring preferential attachment in evolving networks.
\newblock \emph{Europhysics Letters ({EPL})}, 61\penalty0 (4):\penalty0
  567--572, 2003.

\bibitem[Junker and Schreiber(2011)]{junkeranalysis}
Bj{\"o}rn~H Junker and Falk Schreiber.
\newblock \emph{Analysis of biological networks}, volume~2.
\newblock John Wiley \& Sons, 2011.

\bibitem[Kiani et~al.(2015)Kiani, Mahdavi, and Keshavarzi]{kianianalysis}
Rasoul Kiani, Siamak Mahdavi, and Amin Keshavarzi.
\newblock Analysis and prediction of crimes by clustering and classification.
\newblock \emph{Analysis}, 4\penalty0 (8), 2015.

\bibitem[Ku{\v{c}}era(1995)]{kucera}
Lud{\v{e}}k Ku{\v{c}}era.
\newblock Expected complexity of graph partitioning problems.
\newblock \emph{Discrete Applied Mathematics}, 57\penalty0 (2):\penalty0
  193--212, 1995.

\bibitem[Lee and Wilkinson(2019)]{leereview}
Clement Lee and Darren~J Wilkinson.
\newblock A review of stochastic block models and extensions for graph
  clustering.
\newblock \emph{Applied Network Science}, 4\penalty0 (1):\penalty0 122, 2019.

\bibitem[Liu and Wong(2008)]{liu2008effective}
Guimei Liu and Limsoon Wong.
\newblock Effective pruning techniques for mining quasi-cliques.
\newblock \emph{Machine Learning and Knowledge Discovery in Databases}, pages
  33--49, 2008.

\bibitem[Luce and Perry(1949)]{lucemethod}
Robert~D Luce and Albert~D Perry.
\newblock A method of matrix analysis of group structure.
\newblock \emph{Psychometrika}, 14\penalty0 (2):\penalty0 95--116, 1949.

\bibitem[Marinelli et~al.(2020)Marinelli, Pizzuti, and Rossi]{marinellilp}
Fabrizio Marinelli, Andrea Pizzuti, and Fabrizio Rossi.
\newblock Lp-based dual bounds for the maximum quasi-clique problem.
\newblock \emph{Discrete Applied Mathematics}, 2020.

\bibitem[Miao and Balasundaram(2020)]{miaoellipsoidal}
Zhuqi Miao and Balabhaskar Balasundaram.
\newblock An ellipsoidal bounding scheme for the quasi-clique number of a
  graph.
\newblock \emph{INFORMS Journal on Computing}, 2020.

\bibitem[O’Donoghue et~al.(2016)O’Donoghue, Chu, Parikh, and Boyd]{scs}
Brendan O’Donoghue, Eric Chu, Neal Parikh, and Stephen Boyd.
\newblock Conic optimization via operator splitting and homogeneous self-dual
  embedding.
\newblock \emph{Journal of Optimization Theory and Applications}, 169\penalty0
  (3):\penalty0 1042--1068, 2016.

\bibitem[Pardalos and Rodgers(1992)]{pardalos}
Panos~M Pardalos and Gregory~P Rodgers.
\newblock A branch and bound algorithm for the maximum clique problem.
\newblock \emph{Computers \& operations research}, 19\penalty0 (5):\penalty0
  363--375, 1992.

\bibitem[Pardalos and Xue(1994)]{pardalosmaximum}
Panos~M Pardalos and Jue Xue.
\newblock The maximum clique problem.
\newblock \emph{Journal of global Optimization}, 4\penalty0 (3):\penalty0
  301--328, 1994.

\bibitem[Pastukhov et~al.(2018)Pastukhov, Veremyev, Boginski, and
  Prokopyev]{pastukhovmaximum}
Grigory Pastukhov, Alexander Veremyev, Vladimir Boginski, and Oleg~A Prokopyev.
\newblock On maximum degree-based-quasi-clique problem: Complexity and exact
  approaches.
\newblock \emph{Networks}, 71\penalty0 (2):\penalty0 136--152, 2018.

\bibitem[Pattillo et~al.(2012)Pattillo, Youssef, and
  Butenko]{pattillocliquerelaxation}
Jeffrey Pattillo, Nataly Youssef, and Sergiy Butenko.
\newblock Clique relaxation models in social network analysis.
\newblock In \emph{Handbook of Optimization in Complex Networks}, pages
  143--162. Springer, 2012.

\bibitem[Pattillo et~al.(2013)Pattillo, Veremyev, Butenko, and
  Boginski]{pattillomax}
Jeffrey Pattillo, Alexander Veremyev, Sergiy Butenko, and Vladimir Boginski.
\newblock On the maximum quasi-clique problem.
\newblock \emph{Discrete Applied Mathematics}, 161\penalty0 (1):\penalty0
  244--257, 2013.

\bibitem[Pei et~al.(2005)Pei, Jiang, and Zhang]{peimining}
Jian Pei, Daxin Jiang, and Aidong Zhang.
\newblock On mining cross-graph quasi-cliques.
\newblock In \emph{Proceedings of the eleventh ACM SIGKDD international
  conference on Knowledge discovery in data mining}, pages 228--238. ACM, 2005.

\bibitem[Recht et~al.(2010)Recht, Fazel, and Parrilo]{rechtguaranteed}
Benjamin Recht, Maryam Fazel, and Pablo~A Parrilo.
\newblock Guaranteed minimum-rank solutions of linear matrix equations via
  nuclear norm minimization.
\newblock \emph{SIAM review}, 52\penalty0 (3):\penalty0 471--501, 2010.

\bibitem[Ribeiro and Riveaux(2019)]{ribeiroexact}
Celso~C Ribeiro and Jos{\'e}~A Riveaux.
\newblock An exact algorithm for the maximum quasi-clique problem.
\newblock \emph{International Transactions in Operational Research},
  26\penalty0 (6):\penalty0 2199--2229, 2019.

\bibitem[Sturm(1999)]{sturmusing}
Jos~F Sturm.
\newblock Using sedumi 1.02, a matlab toolbox for optimization over symmetric
  cones.
\newblock \emph{Optimization methods and software}, 11\penalty0 (1-4):\penalty0
  625--653, 1999.

\bibitem[Terveen et~al.(1999)Terveen, Hill, and Amento]{terveenconstructing}
Loren Terveen, Will Hill, and Brian Amento.
\newblock Constructing, organizing, and visualizing collections of topically
  related web resources.
\newblock \emph{ACM Transactions on Computer-Human Interaction (TOCHI)},
  6\penalty0 (1):\penalty0 67--94, 1999.

\bibitem[Toh and Yun(2010)]{tohaccelerated}
Kim-Chuan Toh and Sangwoon Yun.
\newblock An accelerated proximal gradient algorithm for nuclear norm
  regularized linear least squares problems.
\newblock \emph{Pacific Journal of optimization}, 6\penalty0
  (615-640):\penalty0 15, 2010.

\bibitem[Toh et~al.(1999)Toh, Todd, and T{\"u}t{\"u}nc{\"u}]{tohsdpt3}
Kim-Chuan Toh, Michael~J Todd, and Reha~H T{\"u}t{\"u}nc{\"u}.
\newblock Sdpt3—a matlab software package for semidefinite programming,
  version 1.3.
\newblock \emph{Optimization methods and software}, 11\penalty0 (1-4):\penalty0
  545--581, 1999.

\bibitem[Veremyev et~al.(2016)Veremyev, Prokopyev, Butenko, and
  Pasiliao]{veremyevexact}
Alexander Veremyev, Oleg~A Prokopyev, Sergiy Butenko, and Eduardo~L Pasiliao.
\newblock Exact mip-based approaches for finding maximum quasi-cliques and
  dense subgraphs.
\newblock \emph{Computational Optimization and Applications}, 64\penalty0
  (1):\penalty0 177--214, 2016.

\bibitem[Verma and Butenko(2013)]{vermanetwork}
Anurag Verma and Sergiy Butenko.
\newblock Network clustering via clique relaxations: A community based.
\newblock \emph{Graph Partitioning and Graph Clustering}, 588:\penalty0 129,
  2013.

\bibitem[Yuan and Yang(2009)]{yuansparse}
Xiaoming Yuan and Junfeng Yang.
\newblock Sparse and low-rank matrix decomposition via alternating direction
  methods.
\newblock \emph{preprint}, 12:\penalty0 2, 2009.

\bibitem[Zhou et~al.(2020)Zhou, Benlic, and Wu]{zhouopposition}
Qing Zhou, Una Benlic, and Qinghua Wu.
\newblock An opposition-based memetic algorithm for the maximum quasi-clique
  problem.
\newblock \emph{European Journal of Operational Research}, 2020.

\end{thebibliography}
\end{document}